\documentclass[12pt,leqno,a4paper]{article}
\usepackage{amsmath,amsthm}
\usepackage[latin1]{inputenc}
\usepackage[T1]{fontenc}
\usepackage[english]{babel}

\newtheorem{theorem}{Theorem}[section]
\newtheorem{vartheorem}{Theorem\kern -2pt}[section]
\newtheorem{proposition}[theorem]{Proposition}
\newtheorem{lemma}[theorem]{Lemma}
\newtheorem{corollary}[theorem]{Corollary}
\newtheorem{definition}[theorem]{Definition}
\newtheorem{conjecture}[vartheorem]{Conjecture\kern -2pt}
\newtheorem{example}[theorem]{Example}
\newtheorem{remark}[theorem]{Remark}
\newtheorem{maintheorem}{Main~Theorem\kern0.5pt}[section]



\newcommand*\varproofname{\textbf{Proof of the main theorem}}
\newcommand*\proofnamestyle{\itshape}

\makeatletter
\newenvironment{varproof}[1][\varproofname]{\par
   \pushQED{\qed}%
   \normalfont \topsep6\p@\@plus6\p@\relax
   \trivlist
   \item[\hskip\labelsep
         \proofnamestyle
     #1\@addpunct{.}]\ignorespaces
}{%
   \popQED\endtrivlist\@endpefalse
}

\newcommand{\dnf}[3]{\frac{d^{#3}#1}{d#2^{#3}}}

\DeclareMathOperator{\tr}{Tr}

\begin{document}

    \title{Trace functions as Laplace transforms}
      \author{Frank Hansen}
      \date{June 28, 2005\\
      {\tiny Latest minor revision August 2, 2005}}
       \maketitle

      \begin{abstract}
      We study trace functions on the form $ t\to\tr f(A+tB) $ where $ f $ is a real
      function defined on the positive half-line, and $ A $ and $ B $ are
      matrices such that $ A $ is positive definite and $ B $ is positive semi-definite.
      If $ f $ is non-negative and operator monotone decreasing, then such a trace function
      can be written as the Laplace transform of a positive measure.
      The question is related to the Bessis-Moussa-Villani conjecture.\\[1ex]
      Key words: Trace functions, BMV-conjecture.
    \end{abstract}

    \section{Introduction}

    The Gibbs density matrix of a system with Hamiltonian $ H $ in equilibrium and
 temperature $ kT=1/\beta $ is given by $ \exp(-\beta H). $ One would like to study
 perturbations  $ H_0+\lambda H_1 $ of an exactly solvable Hamiltonian $ H_0 $ and see how the thermodynamical
 quantities are changed. This question was studied by Bessis, Moussa and Villani in
 \cite{kn:bessis:1975} where it is noted that the Padé approximant to the partition function
 $ Z(\beta)=\tr\exp(-\beta (H_0+\lambda H_1)) $ may be efficiently calculated, if the function
 \[
 \lambda\to\tr\exp(-\beta(H_0+\lambda H_1))
 \]
 is the Laplace transform of a positive measure. The authors then noted that this is indeed true
 for a system of spinless particles with local interactions bounded from below.
 The statement also holds if $ H_0 $ and $ H_1 $ are commuting operators, or if they are just
 $ 2\times 2 $ matrices. These observations led to the formulation of the following conjecture:

  \begin{conjecture}[BMV] Let $ A $ and $ B $ be $ n\times n $ matrices for some natural number
  $ n, $ and suppose that $ A $ is self-adjoint and $ B $ is positive semi-definite. Then there
   is a positive measure $ \mu $  with support in the closed positive half-axis such that
    \[
    \tr\exp(A-tB)=\int_0^\infty e^{-ts}\,d\mu(s)
    \]
    for every $ t\ge 0. $
    \end{conjecture}

    The Bessis-Moussa-Villani (BMV) conjecture may be reformulated as an infinite series of inequalities.

    \begin{vartheorem}[Bernstein] Let $ f $ be a real $ C^\infty $-function defined on the positive half-axis.
    If $ f $ is completely monotone, that is
    \[
    (-1)^n f^{(n)}(t)\ge 0\qquad t>0,\, n=0,1,2,\dots,
    \]
    then there exists a positive measure $ \mu $ on the positive half-axis such that
    \[
    f(t)=\int_0^\infty e^{-st}\, d\mu(s)
    \]
    for every $ t>0. $
    \end{vartheorem}

    The BMV-conjecture is thus equivalent to saying that the function
    \[
    f(t)=\tr\exp(A-tB)\qquad t>0
    \]
    is completely monotone. A proof of Bernstein's theorem can be found in \cite{kn:donoghue:1974}.

    Assuming the BMV-conjecture one may derive a similar statement for free
    semicircularly distributed elements in a type $ II_1 $ von Neumann algebra with a faithful trace. This
    consequence of the conjecture has been proved by Fannes and Petz \cite{kn:fannes:2002}. A hypergeometric
    approach by Drmota, Schachermayer and Teichmann \cite{kn:drmota:2004}
    gives a proof of the BMV-conjecture for some types of $ 3\times 3 $ matrices.

    \subsection{Equivalent formulations}

    The BMV-conjecture can be stated in several equivalent forms.

    \begin{theorem} The following conditions are equivalent:

    \renewcommand{\theenumi}{(\roman{enumi})}

    \begin{enumerate}

    \item For arbitrary $ n\times n $ matrices $ A $ and $ B $
    such that $ A $ is self-adjoint and $ B $ is positive semi-definite the
    function $ f(t)=\tr\exp(A-tB), $ defined on the positive half-axis,
    is the Laplace transform of a positive measure supported in $ [0,\infty). $

    \item For arbitrary $ n\times n $ matrices $ A $ and $ B $
    such that $ A $ is self-adjoint and $ B $ is positive semi-definite the
    function $ g(t)=\tr\exp(A+itB), $ defined on the positive half-axis, is of positive type.

    \item For arbitrary positive definite $ n\times n $ matrices $ A $ and $ B $ the
    polynomial $ P(t)=\tr(A+tB)^p $ has non-negative coefficients for any $ p=1,2,\dots. $

    \item For arbitrary positive definite $ n\times n $ matrices $ A $ and $ B $ the
    function $ \varphi(t)=\tr\exp(A+tB) $ is $ m $-positive on some open interval of the form $ (-\alpha,\alpha). $

    \end{enumerate}

    \end{theorem}

    The first statement is the BMV-conjecture, and it readily implies the second statement by analytic
    continuation. The sufficiency of the second statement is essentially Bochner's theorem. The implication
    $ (iii)\Rightarrow (i) $ is obtained by applying Bernstein's theorem and approximation of the exponential function
    by its Taylor expansion. The implication $ (i)\Rightarrow (iii) $ was proved by Lieb and Seiringer
    \cite{kn:lieb:2003}. A function $ \varphi:(-\alpha,\alpha)\to{\mathbf R} $ is said to be $ m $-positive, if
    for arbitrary self-adjoint $ k\times k $ matrices $ X $ with non-negative entries and spectra contained in
    $ (-\alpha,\alpha) $ the matrix $ \varphi(X) $ has non-negative entries. The implication $ (iii)\Rightarrow (iv) $
    follows by approximation,
    while the implication $ (iv)\Rightarrow (i) $ follows by Bernstein's theorem and
    \cite[Theorem 3.3]{kn:hansen:1992} which states that an $ m $-positive function is real analytic with
    non-negative derivatives in zero.

    \section{Preliminaries and main result}

    Let $ f $ be a real function of one variable defined on a real interval $I.$ We consider for each natural number
    $ n $ the associated matrix function $ x\to f(x) $ defined on the set of self-adjoint matrices of order $ n $
    with spectra in $ I. $ The matrix function is defined by setting
    $$
    f(x)=\sum_{i=1}^p f(\lambda_i) P_i\qquad\mbox{where}\quad
    x=\sum_{i=1}^p \lambda_i P_i
    $$
    is the spectral resolution of $ x. $ The matrix function $ x\to f(x) $
    is Fréchet differentiable~\cite{kn:flett:1980} if $ I $ is open and $ f $ is continuously
    differentiable \cite{kn:brown:1996}.
    The norm of the Fréchet differential $ df(x) $ may be an unbounded function of the order $ m, $
    cf.~\cite{kn:widom:1994, kn:pedersen:2000}. If however $ f $ is assumed to be twice differentiable, then
    the norm of $ df(x) $ is bounded independently of the order $ n $ for all self-adjoint matrices $ x $
    with spectra contained in a fixed compact subset of $ I, $ cf.~\cite[Corollary 2.9]{kn:hansen:1997:2}
    and~\cite{kn:hansen:1995}. We consider in this article the function
    $$
    \varphi(x)=\tr f(x)
    $$
    defined on the set of self-adjoint matrices of order $ n $ with spectra in $ I. $ The Fréchet differential
    is given by $ d\varphi(x)=\tr df(x), $ cf.~\cite{kn:hansen:2002}.

    \subsection{The BMV-property}

    \begin{definition}
    A function $ f\colon \mathbf{R}_+\to\mathbf{R} $ is said to have the BMV-property, if
    to each $ n=1,2,\dots $ and each pair of $ n\times n $
    matrices $ A $ and $ B, $ such that $ A $ is positive definite and $ B $ is positive
    semi-definite, there is a positive measure $ \mu $ with support in $ [0,\infty) $ such that
    \[
    \tr f(A+tB)=\int_0^\infty e^{-st}\, d\mu(s)
    \]
    for every $ t>0. $
    \end{definition}

    The BMV-conjecture is thus equivalent to the statement that the function $ t\to\exp(-t) $ has the
    BMV-property. The main contribution of this paper is the following result.

    \begin{maintheorem}
    Every non-negative operator monotone decreasing function defined on the open positive half-line has the
    BMV-property.
    \end{maintheorem}

    \section{Differential analysis}

    An simple proof of the following result can be found in \cite[Proposition 1.3]{kn:hansen:1995}.

    \begin{proposition}
    The Fréchet differential of the exponential operator function $ x\to\exp(x) $ is given by
    \[
    d\exp(x)h=\int_0^1\exp(sx)h\exp((1-s)x)\,ds=\int_0^1 A(s)\exp(x)\,ds
    \]
    where $ A(s)=\exp(sx)h\exp(-sx) $ for $ s\in{\mathbf R}. $
    \end{proposition}

    This is only a small part of the Dyson formula which contains formalisme developed earlier
    by Tomonaga, Schwinger and Feynman. The subject was given a rigorous mathematical treatment by
    Araki in terms of expansionals in Banach algebras. In particular ~\cite[Theorem 3]{kn:araki:1973},
    the expansional
    \[
    E_r(h; x)=\sum_{n=0}^\infty\int_0^1\int_0^{s_1}\cdots\int_0^{s_{n-1}}A(s_n)A(s_{n-1})\cdots A(s_1)\,
    ds_n\, ds_{n-1}\cdots\, ds_1
    \]
    is absolutely convergent in the norm topology with limit
    \[
    E_r(h; x)=\exp(x+h)\exp(-x).
    \]
    We therefore obtain the $ p $th Fréchet differential of the exponential operator function
    by the expression
    \[
    \begin{array}{l}
    d^p\exp(x)h^p\\[2ex]
    =\displaystyle p!\int_0^1\int_0^{s_1}\cdots\int_0^{s_{p-1}}A(s_p)A(s_{p-1})\cdots A(s_1)\exp(x)
    \,ds_p\,ds_{p-1}\cdots ds_1.
    \end{array}
    \]

    \subsection{Divided differences}

     The following representation of divided differences is due to Hermite
     \cite{kn:hermite:1859}, confer also \cite{kn:norlund:1924, kn:nielsen:1897}.

    \begin{proposition}\label{Hermite expression}
    Divided differences can be written in the following form
    \[
    \begin{array}{rl}
    [x_0, x_1]_f&=\displaystyle\int_0^1 f'\Bigl((1-t_1)x_0+t_1x_1\Bigr)\,dt\\[2ex]
    [x_0, x_1, x_2]_f&=\displaystyle\int_0^1\int_0^{t_1}
    f''\Bigl((1-t_1)x_0+(t_1-t_2)x_1+t_2 x_2\Bigr)\,dt_2\,dt_1\\[2ex]
    &\vdots\\[0pt]
    [x_0, x_1,\cdots, x_n]_f&=\displaystyle\int_0^1\int_0^{t_1}\cdots\int_0^{t_{n-1}}%
    f^{(n)}\Bigl((1-t_1)x_0+(t_1-t_2)x_1+\cdots\\[1ex]
    &\hskip 7em+\,(t_{n-1}-t_n)x_{n-1}+t_nx_n\Bigr)\,dt_n\cdots dt_2\,dt_1
    \end{array}
    \]
    where $ f $ is an $ n $-times continuously differential function defined on an open interval $ I, $ and
    $ x_0,x_1,\dots,x_n $ are (not necessarily distinct) points in $ I. $
    \end{proposition}

    \begin{lemma}\label{lemma: by induction}
    For  $ k=1,2,\dots $ and real numbers $ \lambda_1,\dots,\lambda_k $ we have
    \[
    \int_0^t e^{-\mu s} [s\lambda_1,\dots,s\lambda_k]_{\exp}\, s^{k-1}\,ds
    =t^k e^{-\mu t} [t\lambda_1,\dots,t\lambda_k,t\mu]_{\exp}
    \]
    for any real $ \mu $ and $ t\ge 0. $
    \end{lemma}

    \begin{proof} For $ k=1 $ we calculate
     \[
    \int_0^t e^{-\mu s} [s\lambda_1]_{\exp}\,ds=\int_0^t e^{s(\lambda_1-\mu)}\,ds
    =t\, e^{-\mu t} [t\lambda_1,t\mu]_{\exp}.
    \]
    Assuming the formula valid for $ k $ we obtain for $ k+1 $ the expression
    \[
    \begin{array}{l}
    \displaystyle\int_0^t e^{-\mu s} [s\lambda_1,\dots,s\lambda_{k+1}]_{\exp}\, s^k\,ds\\[2ex]
    =\displaystyle\int_0^t e^{-\mu s}
    \frac{[s\lambda_1,\dots,s\lambda_k]_{\exp}-[s\lambda_2,\dots,s\lambda_{k+1}]_{\exp}}
    {s\lambda_1-s\lambda_{k+1}}\,s^k\,ds\\[2ex]
    =\displaystyle\frac{1}{\lambda_1-\lambda_{k+1}}\Bigl(t^k e^{-\mu t} [t\lambda_1,\dots,t\lambda_k,t\mu]_{\exp}
    -t^k e^{-\mu t} [t\lambda_2,\dots,t\lambda_{k+1},t\mu]_{\exp}\Bigr)\\[3ex]
    =t^{k+1} e^{-\mu t} [t\lambda_1,\dots,t\lambda_{k+1},t\mu]_{\exp}
    \end{array}
    \]
    provided $ \lambda_1\ne\lambda_{k+1}. $ The case $ \lambda_1=\lambda_{k+1} $ then follows by continuity, and
    the lemma is proved by induction.
    \end{proof}

    \begin{theorem}\label{theorem: Frechet differentials of the exponential}
    Let $ x $ and $ h $ be operators on a Hilbert space of finite dimension $ n $
    written on the form
    \[
    x=\sum_{i=1}^n \lambda_i e_{ii}\quad\mbox{and}\quad h=\sum_{i,j=1}^n h_{ij}e_{ij}
    \]
    where $ \{e_{ij}\}_{i,j=1}^n $ is a system of matrix units, and $ \lambda_1,\dots,\lambda_n $ and
    $ h_{i,j} $ for $ i,j=1,\dots,n $ are complex numbers. Then the $ p $th derivative
    \[
    \begin{array}{l}
    \displaystyle\dnf{}{t}{p}\tr \exp({x+th})\Bigr|_{t=0}\\[2ex]
    =\displaystyle p!\sum_{i_1=1}^n\cdots\sum_{i_p=1}^n
    h_{i_p i_{p-1}}\cdots h_{i_2 i_1}h_{i_1 i_p}\,
    [\lambda_{i_1},\lambda_{i_2},\cdots,\lambda_{i_p},\lambda_{i_p}]_{\exp},
    \end{array}
    \]
    where $ [\lambda_{i_1},\lambda_{i_2},\cdots,\lambda_{i_p},\lambda_{i_p}]_{\exp} $ are
    divided differences of order $ p+1 $ of the exponential function.
    \end{theorem}

    \begin{proof}
    We first prove the formulas
    \begin{equation}\label{equations:}
    \begin{array}{rl}
    I_k(s_0)&=\displaystyle\int_0^{s_0}\int_0^{s_1}\cdots\int_0^{s_{k-1}}A(s_k)A(s_{k-1})\cdots A(s_1)
    \,ds_k\,ds_{k-1}\cdots ds_1\\[2ex]
    &=\displaystyle s_0^k \sum_{j=1}^n\sum_{i_1=1}^n\cdots\sum_{i_k=1}^n
    h_{i_k i_{k-1}} h_{i_{k-1} i_{k-2}}\cdots h_{i_2 i_1} h_{i_1 j}\\[1ex]
    &\hskip 10em\times\, e^{-s_0\lambda_j}
    [s_0\lambda_{i_1},\dots,s_0\lambda_{i_k},s_0\lambda_j]_{\exp}\,e_{i_k j}
    \end{array}
    \end{equation}
    for natural numbers $ k=1,\dots,p $\, and real numbers $ s_0\ge 0. $ For $ k=1 $ we calculate the integral
    \[
    \begin{array}{rl}
    \displaystyle I_1(s_0)=\int_0^{s_0} A(s_1)\, ds_1=&\displaystyle\int_0^{s_0}\exp(s_1 x)h\exp(-s_1 x)\,ds_1\\[3ex]
    =&\displaystyle\sum_{i_1=1}^n\sum_{j=1}^n
    h_{i_1 j}\int_0^{s_0}\exp(s_1 x)e_{i_1 j}\exp(-s_1 x)\,ds_1\\[3ex]
    =&\displaystyle\sum_{i_1=1}^n\sum_{j=1}^n
    h_{i_1 j}\int_0^{s_0}\exp(s_1 (\lambda_{i_1}-\lambda_{j})) e_{i_1 j}\,ds_1 \\[3ex]
    =&\displaystyle s_0\sum_{j=1}^n\sum_{1_1=1}^n
    h_{i_1 j}\, e^{-s_0\lambda_j}[s_0\lambda_{i_1}, s_0\lambda_{j}]_{\exp}\, e_{i_1 j}
    \end{array}
    \]
    in accordance with (\ref{equations:}). For $ k\ge 2 $ and assuming the formulas
    $ (\ref{equations:}) $ valid for $ k-1 $ we obtain the expression
    \[
    \begin{array}{rl}
    I_k(s_0)&=\displaystyle\int_0^{s_0}\int_0^{s_1}\cdots\int_0^{s_{k-1}}A(s_k)A(s_{k-1})\cdots A(s_1)
    \,ds_k\cdots ds_2\, ds_1\\[2ex]
    &=\displaystyle\int_0^{s_0} I_{k-1}(s_1) A(s_1)\,ds_1\\[2ex]
    &=\displaystyle\int_0^{s_0} s_1^{k-1}\sum_{m,i_2,\dots,i_k=1}^n
    h_{i_k i_{k-1}}\cdots h_{i_3 i_2} h_{i_2 m}e^{-s_1\lambda_m}\\[1.5ex]
    &\hskip 8em\times\,
    [s_1\lambda_{i_2},\dots,s_1\lambda_{i_k},s_1\lambda_m]_{\exp}\,e_{i_k m} A(s_1)\,ds_1.
    \end{array}
    \]
    We then insert
    \[
    A(s_1)=\displaystyle\sum_{i_1=1}^n\sum_{j=1}^n
    h_{i_1 j}\exp(s_1(\lambda_{i_1}-\lambda_j))\, e_{i_1 j}
    \]
    and using $ e_{i_k m}\, e_{i_1 j}=\delta(m,i_1) e_{i_k j} $ and the symmetry of the divided
    difference we obtain the expression
    \[
    \displaystyle\sum_{j,i_1,i_2,\dots,i_k=1}^n h_{i_k i_{k-1}}\cdots h_{i_2 i_1} h_{i_1 j}
    \int_0^{s_0} e^{-s_1\lambda_j}[s_1\lambda_{i_1},s_1\lambda_{i_2},\dots,s_1\lambda_{i_k}]_{\exp}\,
    s_1^{k-1}\,ds_1\,e_{i_k j}
    \]
    for $ I_k(s_0). $ Finally, using Lemma~\ref{lemma: by induction} we calculate
    \[
    I_k(s_0)=s_0^k\sum_{j,i_1,i_2,\dots,i_k=1}^n h_{i_k i_{k-1}}\cdots h_{i_2 i_1} h_{i_1 j}
     e^{-s_0\lambda_j}[s_0\lambda_{i_1},\dots,s_0\lambda_{i_k},s_0\lambda_j]_{\exp}\,e_{i_k j}
    \]
    which proves (\ref{equations:}) by induction. We next observe that
    \[
    d^p\exp(x)h^p=p!\, I_p(1)\exp(x)
    \]
    where differentiation is with respect to $ x $ when nothing else is indicated. Finally, since
    \[
    \dnf{}{t}{p}\exp({x+th})=d^p_{x+th}\exp(x+th)h^p
    \]
    we obtain
    \[
    \begin{array}{l}
    \displaystyle
    \dnf{}{t}{p}\tr \exp({x+th})\Bigl|_{t=0}=\tr\Big(d^p\exp(x)h^p\Bigr)=p!\,\tr\Bigl(I_p(1)\exp(x)\Bigr)\\[2ex]
    =\displaystyle p!\,\tr\left[\sum_{j,i_1,i_2,\dots,i_p=1}^n h_{i_p i_{p-1}}\cdots h_{i_2 i_1} h_{i_1 j}
     e^{-\lambda_j}[\lambda_{i_1},\dots,\lambda_{i_p},\lambda_j]_{\exp}\,e_{i_p j}\exp(x)\right]\\[3ex]
     =\displaystyle p!\sum_{i_1,i_2,\dots,i_p=1}^n h_{i_p i_{p-1}}\cdots h_{i_2 i_1} h_{i_1 i_p}
     [\lambda_{i_1},\dots,\lambda_{i_p},\lambda_{i_p}]_{\exp}
    \end{array}
    \]
    which is the statement of the theorem.
    \end{proof}

    \begin{lemma}\label{Frechet differential in terms of Fourier tranforms}
    Let $ f\in C_\infty({\mathbf R}), $ and let $ x $ and $ h $ be self-adjoint operators on a
    (possibly infinite dimensional) Hilbert space $ H. $ Then the operator function
    $ x\to f(x) $ is infinitely Fréchet differentiable and the $ p $th Fréchet differential
    is for $ p\ge 1 $ given by
    \[
    d^pf(x)h^p=\int_{-\infty}^\infty d_x^p\exp (-isx)h^p \tilde f(s)\, ds,
    \]
    where $ d_x $ indicates differentiation with respect to $ x $ and
    \[
    \tilde f(s)=\frac{1}{2\pi}\int_{-\infty}^\infty \exp(its) f(t)\, dt
    \]
    is the Fourier transform of $ f. $
    \end{lemma}

    \begin{proof} We note that the statement is true for $ p=1 $ by~\cite[Theorem 1.5]{kn:hansen:1995}
    and assume the statement of the lemma to be valid for $ p. $
    It follows from the definition of the Fréchet differential that the expression
    \[
    d_{x+h}^{p}\exp(-is(x+h))h^p-d_x^{p}\exp(-isx)h^p-d_x(d_x^p\exp(-isx)h^p)h
    \]
    even after division by $ \|h\| $ tend to zero as $ h\to 0. $ We then multiply the above expression
    by the Fourier transform $ \tilde f $ and integrate.
    By Lebesgues's theorem of dominated convergence we therefore obtain that also the expression
    \[
    d_{x+h}^p f(x+h)h^p-d^p_x f(x)h^p-\int_{-\infty}^\infty d_x(d_x^p\exp(-isx)h^p)h\tilde f(s)\, ds,
    \]
    even after division by $ \|h\|, $ tend to zero as $ h\to 0. $ Hence
    \[
    d^{p+1}_x f(x)h^{p+1}=\int_{-\infty}^\infty d_x^{p+1}\exp (-isx)h^{p+1} \tilde f(s)\, ds
    \]
    and the lemma is proved by induction.
    \end{proof}

     In the next corollary we need the identity,
    \begin{gather}\label{identity for divided differences}
    t^{p-1} [t\lambda_1,\dots,t\lambda_p]_f=[\lambda_1,\dots,\lambda_p]_{f_t}
    \qquad\mbox{where}\quad f_t(s)=f(ts),
    \end{gather}
    valid for $ p $ times continuously differentiable functions $ f. $ The statement is
    easily proved by induction after $ p. $

    \begin{corollary}\label{corollary: differentials of trace function}
    Let $ f:I\to{\mathbf R} $ be a $ C^\infty $-function defined on an open and bounded interval $ I, $
    and let $ x $ and $ h $ be self-adjoint operators on a Hilbert space of finite dimension $ n $
    written on the form
    \[
    x=\sum_{i=1}^n \lambda_i e_{ii}\quad\mbox{and}\quad h=\sum_{i,j=1}^n h_{ij}e_{ij}
    \]
    where $ \{e_{ij}\}_{i,j=1}^n $ is a system of matrix units, and $ \lambda_1,\dots,\lambda_n $ are the
    eigenvalues of $ x $ counted with multiplicity. If the spectrum of $ x $ is in $ I, $  then the trace
    function $ t\to\tr f(x+th) $ is infinitely differentiable in a neighborhood of zero and the $ p $th derivative
    \[
    \begin{array}{l}
    \displaystyle\dnf{}{t}{p}\tr f({x+th})\Bigr|_{t=0}\\[2ex]
    \displaystyle=p!\sum_{i_1=1}^n\cdots\sum_{i_p=1}^n
    h_{i_1 i_2}h_{i_2 i_3}\cdots h_{i_{p-1} i_p}h_{i_p i_1}\,
    [\lambda_{i_1},\lambda_{i_2},\cdots,\lambda_{i_p},\lambda_{i_1}]_f,
    \end{array}
    \]
    where $ [\lambda_{i_1},\lambda_{i_2},\cdots,\lambda_{i_p},\lambda_{i_1}]_f $ are
    divided differences of order $ p+1 $ of the function $ f. $
    \end{corollary}

    \begin{proof} Since the spectrum $ \mbox{\textit{Sp}}(x) $ is compact, there is an open and bounded interval
    $ J $ such that
    \[
    \mbox{\textit{Sp}}(x)\subset J\subset\bar{J}\subset I
    \]
    and we may extend the restriction $ f|_J $ to a function $ g\in C_\infty({\mathbf R}). $
    Since the spectrum $ \mbox{\textit{Sp}}(x+th) $ is contained in $ J $ for small $ t $ we
    obtain
    \[
    \begin{array}{l}
    \displaystyle\dnf{}{t}{p}\tr f({x+th})\Bigr|_{t=0}=\dnf{}{t}{p}\tr g({x+th})\Bigr|_{t=0}\\[3ex]
    =\displaystyle\int_{-\infty}^\infty \dnf{}{t}{p}\tr\exp (-isx-isth))\Bigr|_{t=0}\tilde g(s)\, ds\\[3ex]
    =\displaystyle\int_{-\infty}^\infty p!\sum_{i_1=1}^n\cdots\sum_{i_p=1}^n
    (-is)^p h_{i_p i_{p-1}}\cdots h_{i_2 i_1}h_{i_1 i_p}\\
    \hskip 10em\times\,[-is\lambda_{i_1},\cdots,-is\lambda_{i_p},-is\lambda_{i_p}]_{\exp}\,\tilde g(s)\, ds\\[3ex]
    =\displaystyle p!\sum_{i_1=1}^n\cdots\sum_{i_p=1}^n
    h_{i_p i_{p-1}}\cdots h_{i_2 i_1}h_{i_1 i_p}\int_{-\infty}^\infty
    [\lambda_{i_1},\cdots,\lambda_{i_p},\lambda_{i_p}]_{\exp(-is\cdot)}\,\tilde g(s)\, ds\\[3ex]
    =\displaystyle p!\sum_{i_1=1}^n\cdots\sum_{i_p=1}^n
    h_{i_p i_{p-1}}\cdots h_{i_2 i_1}h_{i_1 i_p}\,
    [\lambda_{i_1},\cdots,\lambda_{i_p},\lambda_{i_p}]_g\, ds,
    \end{array}
    \]
    where $ \tilde g $ is the Fourier transform of $ g, $ and we used
    (\ref{identity for divided differences}) and the linearity in $ g $ of arbitrary divided
    differences $ [\lambda_1,\dots,\lambda_p]_g. $ The statement now follows by noting that
    $ f=g $ in a neighborhood of the spectrum of $ x. $
    \end{proof}

    \section{Proof of the main theorem}

    \begin{proposition}
    Consider for a constant $ c\ge 0 $ the function
    \[
    g(t)=\frac{1}{c+t}\qquad t>0.
    \]
    For arbitrary $ n\times n $ matrices $ x $ and $ h $ such that
    $ x $ is positive definite and $ h $ is positive semi-definite
    we have
    \[
    (-1)^p  \dnf{}{t}{p}\tr g({x+th})\Bigr|_{t=0}\ge 0
    \]
    for $ p=1,2,\dots. $
    \end{proposition}

    \begin{proof} We first note that the divided differences of $ g $ are of the form
    \begin{gather}
    [\lambda_1,\lambda_2,\dots,\lambda_p]_g=(-1)^{p-1}g(\lambda_1)g(\lambda_2)\cdots g(\lambda_p)\qquad
    p=1,2,\dots
    \end{gather}
    for arbitrary positive numbers $ \lambda_1,\lambda_2,\dots,\lambda_p. $
    There is nothing to prove for $ p=1. $ Assume the statement true for $ p\ge 2 $ and notice
    that
    \[
    \frac{g(\lambda)-g(\mu)}{\lambda-\mu}=\frac{(c+\lambda)^{-1}-(c+\mu)^{-1}}{\lambda-\mu}=-g(\lambda)g(\mu)
    \]
    for $ \lambda\ne\mu. $ Therefore
    \[
    \begin{array}{rl}
     [\lambda_1,\dots,\lambda_p,\lambda_{p+1}]_g
     &=\displaystyle\frac{[\lambda_1,\dots,\lambda_p]_g
    -[\lambda_2,\dots,\lambda_{p+1}]_g}{\lambda_1-\lambda_{p+1}}\\[3ex]
     &=\displaystyle (-1)^{p-1}
     g(\lambda_2)\cdots g(\lambda_p)\frac{g(\lambda_1)-g(\lambda_{p+1})}{\lambda_1-\lambda_{p+1}}\\[3ex]
     &=(-1)^p g(\lambda_1)\cdots g(\lambda_p)g(\lambda_{p+1})
     \end{array}
     \]
     for $ \lambda_1\ne\lambda_{p+1} $ and the general case follows by approximation. Next, applying
     Corollary~\ref{corollary: differentials of trace function} we obtain
    \[
    \begin{array}{l}
    \displaystyle\dnf{}{t}{p}\tr g({x+th})\Bigr|_{t=0}\\[2ex]
    \displaystyle=p!\sum_{i_1=1}^n\cdots\sum_{i_p=1}^n
    h_{i_1 i_2}h_{i_2 i_3}\cdots h_{i_{p-1} i_p}h_{i_p i_1}\,
    [\lambda_{i_1},\lambda_{i_2},\cdots,\lambda_{i_p},\lambda_{i_1}]_g\\[2ex]
    =\displaystyle(-1)^p p!\sum_{i_1=1}^n\cdots\sum_{i_p=1}^n
    h_{i_1 i_2}h_{i_2 i_3}\cdots h_{i_{p-1} i_p}h_{i_p i_1}\,
    g(\lambda_{i_1})g(\lambda_{i_2})\cdots g(\lambda_{i_p})g(\lambda_{i_1}).
    \end{array}
    \]
    Since $ h $ is positive definite, it is of the form $ h=aa^* $ for some matrix $ a $ and therefore
    \[
    h_{ij}=\sum_{m=1}^n a_{im}a^*_{mj}=\sum_{m=1}^n a_{im} \overline{a}_{jm}=(a_i\mid a_j)\qquad i,j=1,\dots,n,
    \]
    where
    \[
    a_i=\left(\begin{array}{ccc}
               a_{i1} & \cdots & a_{in}
             \end{array}\right)\in{\mathbf C}^n
    \]
    is the $ i $th row in the matrix $ a. $
    We then set $ \xi_i=g(\lambda_i) a_i $ and $ b_i=g(\lambda_i)^{1/2} a_i $ for $ i=1,\dots,n $ and calculate
    \[
    \begin{array}{l}
    \displaystyle\frac{(-1)^p}{p!}\,\dnf{}{t}{p}\tr g({x+th})\Bigr|_{t=0}\\[2ex]
    =\displaystyle\sum_{i_1=1}^n\cdots\sum_{i_p=1}^n
    h_{i_1 i_2}h_{i_2 i_3}\cdots h_{i_{p-1} i_p}h_{i_p i_1}\,
    g(\lambda_{i_1})g(\lambda_{i_2})\cdots g(\lambda_{i_p})g(\lambda_{i_1})\\[2ex]
    =\displaystyle\sum_{i_1=1}^n\cdots\sum_{i_p=1}^n
    (\xi_{i_1}\mid b_{i_2})(b_{i_2}\mid b_{i_3})\cdots (b_{i_{p-1}}\mid b_{i_p})(b_{i_p}\mid \xi_{i_1}).
    \end{array}
    \]
    But any sum of the form
    \[
    S=\sum_{i_2=1}^n\cdots\sum_{i_p=1}^n (\xi\mid b_{i_2})(b_{i_2}\mid b_{i_3})
    \cdots (b_{i_{p-1}}\mid b_{i_p})(b_{i_p}\mid \xi)
    \]
    is non-negative. It is obvious for $ p=2 $ since then
    \[
    S=\sum_{i_2=1}^n (\xi\mid b_{i_2})(b_{i_2}\mid \xi),
    \]
    and for $ p=3 $ since then
     \[
    \sum_{i_2=1}^n\sum_{i_3=1}^n (\xi\mid b_{i_2})(b_{i_2}\mid b_{i_3})(b_{i_3}\mid \xi)
    =\left(\sum_{i_2=1}^n (\xi\mid b_{i_2})b_{i_2}\mid \sum_{i_3=1}^n (\xi\mid b_{i_3})b_{i_3}\right).
    \]
    For $ p\ge 4 $ we may therefore use induction over either the even or the odd natural numbers by noting that
    \[
    \begin{array}{rl}
    S&=\displaystyle\sum_{i_2=1}^n\cdots\sum_{i_p=1}^n (\xi\mid b_{i_2})(b_{i_2}\mid b_{i_3})
    \cdots (b_{i_{p-1}}\mid b_{i_p})(b_{i_p}\mid \xi)\\[2ex]
    &=\displaystyle\sum_{i_3=1}^n\cdots\sum_{i_{p-1}=1}^n (\eta\mid b_{i_3})(b_{i_3}\mid b_{i_4})
    \cdots (b_{i_{p-2}}\mid b_{i_{p-1}})(b_{i_{p-1}}\mid \eta)
    \end{array}
    \]
    where the vector
    \[
    \eta=\sum_{i=1}^n (\xi\mid b_i)b_i.
    \]
    This concludes the proof of the statement.
    \end{proof}

    \begin{varproof}
    Consider again the function
    \[
    g(t)=\frac{1}{c+t}\qquad t>0
    \]
    for $ c\ge 0 $
    and arbitrary $ n\times n $ matrices $ x $ and $ h $ such that $ x $
    is positive definite and $ h $ is positive semi-definite. We first note that
    \[
    \dnf{}{t}{p}\tr g({x+th})\Bigr|_{t=t_0}
    =\dnf{}{\varepsilon}{p}\tr g({x+t_0h+\varepsilon h})\Bigr|_{\varepsilon=0}
    \]
    for $ p=1,2,\dots $ and $ t_0\ge 0. $ The function
    $ t\to\tr g(x+th) $ is therefore completely monotone.
    Let now $ f\colon{\mathbf R}_+\to\mathbf R $ be a non-negative operator monotone decreasing function.
    Any operator monotone decreasing function defined on the open positive half-line, thus in particular
    the function $ f, $ is necessarily of the form
    \begin{gather}\label{operator monotone decreasing function}
    f(t)=\alpha t +\beta + \int_0^\infty\left(\frac{1}{c+t}-\frac{c}{c^2+1}\right)\,d\nu(c),
    \end{gather}
    where $ \alpha\le 0 $ and $ \nu $ is some positive Borel measure with support in $ [0,\infty) $ such that the
    integral $ \int (c^2+1)^{-1}\,d\nu(c) $ is finite,
    cf.~\cite[Chap.~II Theorem~1 and Lemma~2]{kn:donoghue:1974}. Note that we may write
    \[
    \int_0^\infty\left(\frac{1}{c+t}-\frac{c}{c^2+1}\right)\,d\nu(c)
    =\int_0^\infty\frac{1-ct}{c+t}\cdot(c^2+1)^{-1}\,d\nu(c)
    \]
    where for each $ t>0 $ the function $ c\to (1-ct)(c+t)^{-1} $ is decreasing and
    bounded between $ -t $ and $ t^{-1}. $
    Since $ f $ has a finite limit
    for $ t\to\infty $ one may derive that $ \alpha=0, $ and by appealing to Lebesgues's theorem
    of monotone convergence we also obtain $ \int c(c^2+1)^{-1}\,d\nu(c)<\infty, $
    hence $ f $ allows the representation
    \[
    f(t)=\beta+\int_0^\infty\frac{1}{c+t}\,d\nu(c),
    \]
    where we have incorporated the constant contribution from the integral
    in (\ref{operator monotone decreasing function}) into $ \beta $ such that
    \[
    \beta=\lim_{t\to\infty} f(t)\ge 0.
    \]
    We conclude that the function $ t\to\tr f(x+th) $ is completely monotone
    and thus by Bernstein's theorem is the Laplace transform of a positive measure with support in $ [0,\infty). $
    \end{varproof}

    \subsection{Further analysis}

    One may try to use the Hermite expression in Proposition~\ref{Hermite expression} to obtain a proof of the
    BMV-conjecture. Applying Theorem~\ref{theorem: Frechet differentials of the exponential} and calculating
    the third derivative of the trace function we obtain
    \[
    \begin{array}{l}
    \displaystyle\frac{-1}{3!}\,\dnf{}{t}{3}\tr\exp({x-th})\Bigr|_{t=0}
    \displaystyle=\sum_{p,i,j=1}^n
    (a_p\mid a_i) (a_i\mid a_j) (a_j\mid a_p)[\lambda_p\lambda_i\lambda_j\lambda_p]_{\exp}\\[3ex]
    =\displaystyle\int_0^1\int_0^{t_1}\int_0^{t_2}\sum_{p,i,j=1}^n
    {(a_p\mid a_i) (a_i\mid a_j) (a_j\mid a_p)}
    \displaystyle\exp\bigl((1-(t_1-t_3))\lambda_p\\
    \hskip 15em+\,(t_1-t_2)\lambda_i+(t_2-t_3)\lambda_j\bigr)\,dt_3\,dt_2\,dt_1
    \end{array}
    \]
    where $ h=aa^* $ and $ a_i $ is the $ i $th row in $ a. $ Assuming the BMV-conjecture this
    integral should be non-negative, and this would obviously be the case if the integrand is a
    non-negative function.

    \begin{example}\rm If we evaluate the above integrand in $ t_1=t_2=1 $ and $ t_3=1/3 $ we obtain
    \[
    S=\sum_{p,i,j=1}^n
    (a_p\mid a_i) (a_i\mid a_j) (a_j\mid a_p)\\
    \exp\Bigl(\frac{\lambda_p+2\lambda_j}{3}\Bigr)\,dt_3\,dt_2\,dt_1.
    \]
    If we in addition set
    \[
    \begin{array}{rl}
    a_1&=(1000, -10, 1)\\[0.5ex]
    a_2&=(-10, 10000, 1000)\\[0.5ex]
    a_3&=(1, 1000, 202139)
    \end{array}
    \]
    and $ (\lambda_1, \lambda_2, \lambda_3)=3\log 2\cdot(23,11,0), $ then the sum is an integer with value
    \[
    S=-487062506352658941731358505750.
    \]
    The values of $ S $ are extremely sensitive to the chosen figures, and they tend to be overwhelmingly
    positive. If for example the third entry in $ a_3 $ is changed from 202139 to 202138, then
    \[
    S=376189230591238013538921396773.
    \]
    The result is equally sensitive to changes in the values of $ (\lambda_1, \lambda_2, \lambda_3). $
    \end{example}

    The above
    example has not been found by simulation. In fact, millions of simulations with randomly chosen entries
    have been carried out, without ever hitting a negative value of $ S. $ Instead,
    the example has been obtained by the study and proper modification of an
    example in~\cite{kn:johnson:2002} of two positive definite $ 3\times 3 $ matrices $ A $ and $ B $ such that
    $ \tr(BABAAB)=-3164. $

    Another way forward would be to examine the value of loops of the form
    \[
    (a_1\mid a_2)(a_2\mid a_3)\cdots(a_{p-1}\mid a_p)(a_p\mid a_1)
    \]
    since they, apart from an alternating sign, are the only possible negative factors in the expression of the derivatives of the
    trace functions. The value of a loop is a homogeneous function of degree two in the norm of the vectors, so we
    only need to consider unit vectors. The value of such a loop is furthermore invariant under unitary
    transformations and thus takes the minimal value in any subspace of dimension $ p. $ By applying a variational
    principle the lower bound
     \[
    -\cos^n\left(\frac{\pi}{n}\right)\le(a_1\mid a_2)(a_2\mid a_3)\cdots(a_{p-1}\mid a_p)(a_p\mid a_1)
    \]
    was established in \cite{kn:hansen:2004:2}. The lower bound converges very slowly to $ -1 $ as $ p $
    tends to infinity, and it is attained essentially only when all
    the vectors form a ''{\it fan}'' in a two-dimensional subspace.

    \begin{remark}
    If we only consider one-dimensional perturbations, that is if $ h=cP $
    for a constant $ c>0 $ and a one-dimensional projection $ P, $
    then $ h $ is of the form $ h=(\xi_i\bar\xi_j)_{i,j=1,\dots,n} $ for some vector
    $ \xi=(\xi_1,\dots,\xi_n) $ and each loop
    \[
    h_{i_1 i_2}h_{i_2 i_3}\cdots h_{i_{p-1} i_p}h_{i_p i_1}
    =\|\xi_{i_1}\|^2\cdots\|\xi_{i_p}\|^2
    \]
    is manifestly real and non-negative. This implies that the trace function
    \[
    t\to\tr\exp(-(x+th)),
    \]
    for any self-adjoint $ n\times n $ matrix $ x, $ is the
    Laplace transform of a positive measure with support in
    $ [0,\infty). $ The same statement holds, with $ x $ positive definite, for the trace function
    $ t\to\tr f(x+th) $ associated with an arbitrary completely monotone function $ f. $

    \end{remark}


      \vfill

      {\small\noindent Frank Hansen: Institute of Economics, University
       of Copenhagen, Studiestraede 6, DK-1455 Copenhagen K, Denmark.}

      \end{document}